\documentclass[12pt]{article}

\usepackage{harvard}
\usepackage{amssymb}
\usepackage{amsmath}

\newtheorem{theorem}{Theorem}
\newtheorem{lemma}{Lemma}
\newtheorem{corollary}{Corollary}
\newtheorem{proposition}{Proposition}
\newtheorem{remark}{Remark}

\newtheorem{example}{Example}

\newcommand{\ignore}[1]{}

\newcommand{\ch}[1]{\citeasnoun{#1}}

\newcommand{\AO}{A_0}
\newcommand{\Al}{A_1}
\newcommand{\Sgmall}{\mathbf{\Sigma}_{11}}
\newcommand{\SgmaOO}{\mathbf{\Sigma}_{00}}
\newcommand{\All}{A_{11}}
\newcommand{\AOO}{A_{00}}
\newcommand{\AllI}{A_{11}^{-1}}
\newcommand{\AOOI}{A_{00}^{-1}}
\newcommand{\Vll}{V_{11}}
\newcommand{\VOO}{V_{00}}
\newcommand{\Yl}{\mathbf{Y}_1}
\newcommand{\YO}{\mathbf{Y}_0}

\newenvironment{ofproof}[3]{{\bf Proof of #1  #2.} \em  \hfill \ \linebreak  \indent #3 $\diamondsuit$ }{}
\renewcommand{\cleardoublepage}
    {\clearpage\if@twoside \ifodd\c@page\else
    \hbox{}\thispagestyle{empty}\newpage\if@twocolumn\hbox{}\newpage\fi\fi\fi}

\begin{document}
\pagestyle{myheadings}
\markright{\underline{VRP for OLS}}

\vspace{-3.0in}

\title{ Does sequential augmenting\\
of the simple linear heteroscedastic regression\\
reduce variances of the Ordinary Least Squares?\footnote{The final and revised version of the paper has been accepted for publication in {\it Statistics, ISSN: 0233-1888 E-ISSN: 1029-4910}}
}

\author{\normalsize Andrzej S. Kozek$^\ast$$^a$\thanks{$^\ast$Corresponding author. Email: Andrzej.Kozek@mq.edu.au
\vspace{6pt}} and Brian Jersky$^b$\\\vspace{6pt}  {\small\em{Macquarie University$^a$, Sydney and California State Polytechnic University$^b$, Pomona}}}
\date{\small 2012-03-22} 
\maketitle
\thispagestyle{empty}
\hrule

\begin{abstract}
 If uncorrelated random variables have a common expected value and decreasing variances then the variance of a sample mean  is decreasing with the number of observations. Unfortunately, this natural and desirable Variance Reduction Property (VRP) by augmenting data is not automatically inherited by Ordinary Least Squares (OLS) estimators of parameters.
 In the paper we find conditions for the OLS to have the VRP. In the case of a straight line regression we show that the OLS estimators of intercept and slope  have the VRP if the design points are increasing. This also holds true for alternating two-point experimental designs. The obtained results are useful in the cases where it is known that variances of the subsequent observations are non-increasing, but the ratios of the decrease are not available to use sub-optimal or optimal Weighted Least Squares estimators. \\\nopagebreak
  {\em AMS (2010) subject classification.} Primary 62J05; secondary 62K99,  15A45, 62G20.
\end{abstract}

\hrule

\section{Introduction} \label{s1}
If uncorrelated random variables have a common expected value and decreasing variances then, in agreement with common sense, the variance of a sample mean  is decreasing with the number of observations. Surprisingly, this natural and desirable Variance Reduction Property (VRP) resulting from inclusion of less contaminated data is not shared in general by the Ordinary Least Squares (OLS) estimators of parameters. 

In fact, by generating random regression designs, it is easy to come across experimental designs where by adding additional observation with a smaller error variance, we in fact increase the  variance of the augmented OLS estimators (we refer to Example \ref{ex03} for a specific case). On the other hand, many such random designs indeed enjoy the VRP. In the case of homoscedastic models with identical variances, or when the covariance matrix is known up to a constant, the VRP is well known,  with the Weighted Least Squares replacing the OLS in the latter case, and a fairly complete theory on optimal augmenting of the data is then available. We refer the reader to \ch{MR1392638} and \ch{MR1897772} for the most comprehensive recent results in this direction, and for further references to the related literature. The case of unknown and decreasing variances appears difficult and seemingly has eluded the attention of researchers. This case however covers important applications in the theory of Asymptotic Statistics, where the asymptotic distributions of smoothed functionals are exactly of such a form.

In the paper we derive (Theorem \ref{generalTh}) general conditions characterizing designs for which the OLS estimators of parameters
retain the variance reduction property when new data with lower variances are included. These criteria can be easily used to iteratively construct sequences of experimental designs with the VRP. We found one interesting general case where the VRP remains valid: when $n$ uncorrelated observations have equal variances of errors and the $(n+1)-$th observation has a lower variance.  Then the variance of the OLS estimator of parameters based on $(n+1)$ observations is lower than that based on $n$ observations (Theorem \ref{W114I}).  Following the method of Theorem \ref{generalTh} in a particularly important linear model with two parameters: the intercept and slope, we show then that  any design with the explanatory variable increasing with $n$ has the VRP property (Theorem \ref{mainTh}), and that any two-point alternating experimental design has the VRP property (Theorem \ref{twopointsdesignth}). All these results are non-trivial, novel and guarantee non-negativeness of the diagonal elements of matrices which themselves need not be positively definite.

The paper is organized in the following way. In Section \ref{general} we present the main results of the paper. Section \ref{linear} contains details of the straight line regression model and proofs of Theorems \ref{mainTh} and \ref{twopointsdesignth}. In the proofs we include only the main nontrivial steps and omit the detailed, often long, though elementary algebraic manipulations. Technical details can be found in  Appendix \ref{app}.
In Section \ref{examples} we include examples illustrating that the assumptions of Theorem \ref{mainTh} cannot be weakened.  
In Appendix \ref{app} technical auxiliary results needed in the proofs of the main theorems of the paper are presented. We also include Proposition \ref{t2dep} showing decomposition of variances in the general case of dependent random variables. Proposition \ref{t2dep} extends the decomposition \eqref{kwadratyporedukcjidoI} valid for uncorrelated random variables. We acknowledge using {\sc Scientific Workplace} for deriving the results of this paper.\\

\section{Variance reduction in heteroscedastic models}\label{general}
Consider the linear model
\begin{align}
\YO  &=\AO \beta +\epsilon_0, \label{m2O}
\end{align}
where $\AO $ is an $(n,k)$ design matrix, $\beta$ is a $(k,1)$ column vector of regression parameters, and $\YO = \left(Y_1, Y_2, \ldots, Y_n \right)^T$ and $\epsilon_0 \in R^{n}$ are random vectors. We also consider an augmented model
\begin{align}
\Yl  &=\Al \beta +\epsilon_1, \label{m2a}
\intertext{where $\Al =\left[
\begin{array}{c}
\AO  \\
a%
\end{array}
\right] $ is an $(n+1,k)$ design matrix with  $a$ being a design of the $(n+1)$-th experiment, $\Yl = \left[\YO^T, Y_{n+1} \right]^T$ and  $\epsilon = \left[ \epsilon_0^T , \epsilon_{n+1}\right]^T \in R^{n+1}$ are random vectors such that}
E\left(\epsilon_1\right)=0,\quad &\quad Cov(\epsilon_1 ) =\Sgmall.  \label{m2b}
\end{align}
Apart from Example \ref{ex02} and Proposition \ref{t2dep} we will consider uncorrelated random variables with a diagonal covariance matrix $\Sgmall$.  Hence we have  
\begin{align}
\Sgmall  &=\left[
\begin{array}{cc}
\SgmaOO & O \\
O^{T } & \sigma_{n+1}^{2}%
\end{array}%
\right], \label{m2c}
\end{align}
where $O$ is an $(n,1)$ vector of zeros and $\SgmaOO$ is a diagonal covariance matrix of $\YO$. To simplify notation we define
\begin{align}
\All &=  \Al ^{T}\Al \quad \text{and} \label{def4All}\\
\AOO &=  \AO ^{T}\AO, \label{def4AOO}
\end{align}
respectively.
Both $\beta$ and the covariance matrix $\Sgmall$ are assumed to be unknown and the design matrix $\Al $ is assumed to be of full rank. The OLS  estimators of parameters $\beta$ based on models \eqref{m2O} and \eqref{m2a} are of the form
\begin{align}
 \hat{\beta_1} &= \AllI \Al ^T  \Yl  \label{def4beta}
 \intertext{and}
 \hat{\beta_0} &= \AOOI \AO ^T  \YO ,  \label{def4beta}
 \intertext{and have covariance matrices $\Vll $ and $\VOO$ given by}
\Vll &= \AllI \Al ^{T}\mathbf{\Sgmall }\Al \AllI \label{def4V00}
\intertext{and}
\VOO &= \AOOI \AO ^T \mathbf{\SgmaOO }\AO \AOOI , \label{def4VOO}
\end{align}
\nopagebreak respectively. In the paper we derive conditions implying inequalities
\begin{equation} \label{varred}
\Vll (i,i) \le \VOO (i,i) \quad \text{for $i\in\{1,2,\ldots,k\},$}
\end{equation}
where $\Vll $ and $\VOO $ are given by \eqref{def4V00}   and \eqref{def4VOO}, respectively. We will refer to property \eqref{varred} as  the variance reduction by augmenting the data, or, for short, as the VRP property.

It will be convenient to use  the following notation throughout the paper
\begin{align}
d_{n+1} &= a \AOOI a^{T}, \quad {and}\label{def4d}\\
q_{n+1} &= 1+d_{n+1}.\label{def4q}
\end{align}

In the general case of a heteroscedastic model we have the following decomposition of variance $\Vll $.
\begin{proposition}
\label{t2} If the covariance matrix $\mathbf{\Sgmall }$ and design matrix $\Al $ in the model given by \eqref{m2a}  are both of full rank,   then for uncorrelated random variables $Y_i$, $i=1,\ldots,n+1$, we have
\begin{align}
\Vll  =
\VOO -\mathbf{\sigma }_{n+1}^{2}W-W_{11}, \label{kwadratyporedukcjidoI}
\end{align}
where
\begin{align}
W &= \frac{1}{q_{n+1}}\AOOI a^{T}a\AOOI ,\label{def4W}\\
W_{11} &=  \AOOI D_{11} \AOOI -\AllI D_{11} \AllI , \label{def4W11}\\
D_{11} &= \AO ^{T}D\AO ,  \quad \text{and} \label{def4D11}\\
D &= \mathbf{\SgmaOO}-\mathbf{\sigma }_{n+1}^{2}I.  \label{def4D}
\end{align}
\end{proposition}

 In the particular case of iid random variables with common variance $\sigma^2$ the matrices $D$ and $W_{11}$ equal zero and then \eqref{kwadratyporedukcjidoI} reduces to the well known case of the \ch{MR0036980} updating formula for covariances of the LS estimators
\begin{equation}
\AllI = \AOOI -W . \label{iidcase}
\end{equation}
Let us note that in \eqref{kwadratyporedukcjidoI}  the matrix $\sigma_{n+1}^2W$ is non-negative definite and hence it always reduces the variances on the diagonal of $\Vll $.  However, in general, the matrix $W_{11}$ need not be positive definite, cf. Example \ref{ex01} in Section \ref{examples}, and for general designs it may even have negative elements on the diagonal. Consequently, in the heteroscedastic case the variance of estimators of $\beta$ can increase with an increase of the number of observations $n$, even  when $k$ and the parameters $\beta$ are kept the same.
As might be expected, if $\sigma_{n+1} > \max\left\{\sigma_i, i=1,\ldots,n \right\}$, then the variance of $\hat{\beta}_1$ could be greater than that of $\hat{\beta}_0$. Therefore in the following we shall concentrate on the cases where  the covariance matrices
\begin{equation}
\Sgmall = Diag\left\{\sigma_{1}^2, \sigma_{2}^2,  \cdots \sigma_{n+1}^2
\right\}\label{diagSigma}
\end{equation}
 have  non-increasing standard deviations $\sigma_{i}, i = 1,2,\ldots,n+1$. Clearly, in these cases the matrix $D$ given by \eqref{def4D} is also diagonal with non-increasing diagonal elements
 $D_i = \sigma^2_i -\sigma^2_{n+1}\ge~0. $
 As has been mentioned previously in the Introduction, decreasing $\sigma_i$'s alone do not guarantee the VRP, however.

 Nonetheless, we have the following surprising general case where the VRP is valid: if the first $n$ observations have uncorrelated errors with the same variances and the error in the $(n+1)$-th uncorrelated experiment is lower than in the preceding cases, then the ordinary least squares estimators based on the augmented data have lower variances than in the case of $n$ observations and the homoscedastic model.
\begin{theorem}\label{W114I}
If  $D=\sigma^2I$ where I is the identity matrix then
\begin{equation} \label{surpriseVRP}
W_{11}=  \sigma^2\left( 2-\frac{d_{n+1}}{1+d_{n+1}}\right) W
\end{equation}
and, consequently, $W_{11}$ is non-negative definite.~$\lozenge$
\end{theorem}

Clearly Theorem \ref{W114I} and Proposition \ref{t2} imply the following corollary.
\begin{corollary} \label{surprise}
If random variables $Y_1,Y_2,\ldots,Y_{n+1}$ are uncorrelated, $var(Y_i)=\sigma_1^2$ for $i=1,2,\ldots,n$ and $var(Y_{n+1})=\sigma_{n+1}^2$ with $\sigma_1^2 >\sigma_{n+1}^2$ then the OLS estimators have the variance reduction property.\ $\lozenge$
\end{corollary}
The proof of Theorem \ref{W114I} can be obtained using \eqref{prosto1} and \eqref{prosto2} and some straightforward 
algebra, and is omitted.

In the general uncorrelated case we have the following characterization of the non-negative diagonal elements of  $W_{11}$, implying the VRP.
\begin{proposition}\label{char4WiiPos}
The diagonal element $W_{11}(i,i)$ of $W_{11}$ is non-negative for every diagonal covariance matrix \eqref{diagSigma} with $\sigma_{i}$ non-increasing if and only if
\begin{equation} \label{techCond4Wiipos}
\sum_{j=1}^m \left(\left(\sum_{l=1}^k \AOOI (i,l)\AO (j,l)\right)^2 -\left(\sum_{l=1}^k \AllI (i,l)\AO (j,l)\right)^2\right) \ge 0 \quad \text{for\ } m=1,2,\ldots,n.
\end{equation}
\end{proposition}

The following theorem follows immediately from Propositions \ref{t2} and \ref{char4WiiPos}.
\begin{theorem}\label{generalTh}
The OLS estimator of the $i-th$ component of $\beta$ has the VRP if condition \eqref{techCond4Wiipos} is met.
\end{theorem}

 Conditions \eqref{techCond4Wiipos} are in general not easy to verify. However, given the design matrix $\Al $ of the first $n$ experiments one can use \eqref{iidcase} to find numerically designs $a$ for the  $(n+1)$-st experiment satisfying inequalities \eqref{techCond4Wiipos}. Hence Proposition \ref{char4WiiPos} provides a method for a sequential numerical construction of experimental designs having the VRP. Remark \ref{r2} in Section \ref{linear} provides more details in the case of a straight line linear model.

 The linear model for a straight line regression has the following design matrix for $n>1$.
\begin{align}
\Al &=\left[ \begin{array}{cc}
  \mathbf{1} & \mathbf{h}_{n+1} \\
\end{array}\right]\label{def4h}
= \left[%
\begin{array}{cc}
  1 & h_1 \\
  \vdots & \vdots \\
  1 & h_n \\
  1 & h_{n+1} \\
\end{array}%
\right].
\end{align}

The following two theorems show that in this important case the VRP remains a valid feature of two important classes of regression designs: with non-negative and increasing explanatory variables $h_i$ (Theorem \ref{mainTh}) and with two alternating values of the explanatory variable (Theorem \ref{twopointsdesignth}).

\begin{theorem}\label{mainTh}
If random variables $Y_i, i = 1,2,\ldots,n+1$ are uncorrelated with non-increasing variances $\sigma_i^2$ and with increasing values of the explanatory variable $h_i\ge 0$
then the OLS estimators with design matrix \eqref{def4h} have the variance reduction property \eqref{varred}.
\end{theorem}

\begin{remark} \label{monoton} In Theorem \ref{mainTh} we consider models with a simple explanatory variable $h$. Clearly, any increasing function of $h$ can be used here, as such a model, by applying a  suitable change of parametrization, can be easily transformed to the case with a design matrix \eqref{def4h}. 
\end{remark}

Examples \ref{ex01} -- \ref{ex03} of Section \ref{examples} show that the assumptions of Theorem \ref{mainTh} cannot be weakened.  This however does not contradict the fact that the VRP may also hold for many regression designs which do not meet the monotonicity assumption of Theorem \ref{mainTh}. For example, it is well known that in the homoscedastic case two-point simple regression designs result in the lowest variance of the LS estimators of the slope and intercept. The following theorem shows that these two-point regression designs also have the variance reduction property.
\begin{theorem}\label{twopointsdesignth}
If random variables $Y_i, i = 1,2,\ldots,n+1$ are uncorrelated with variances $\sigma_i^2$ non-increasing with $i$, with two distinct values of the explanatory variable given by
\begin{equation}\label{2pdes}
h_i=\left\{
      \begin{array}{ccc}
        b & \text{if} & i \text{\ is odd} \nonumber \\
        c & \text{if} & i \text{\ is even} \nonumber  \\
      \end{array}
    \right.
\end{equation}
then the OLS estimators \eqref{def4beta} with design matrix \eqref{def4h} and non-increasing variances \eqref{diagSigma} have the variance reduction property \eqref{varred}.
\end{theorem}

\section{Variance reduction for a straight line model}\label{linear}
Consider the case of the straight line model \eqref{def4h}. 
By Theorem \ref{generalTh} to prove the variance reduction property  we need to verify that conditions \eqref{techCond4Wiipos} hold true. In this particular case \eqref{techCond4Wiipos} is equivalent to  inequalities \eqref{mainIneq1a} and \eqref{mainIneq1a} of Condition C3. In Theorem \ref{equivCond} we show implication relations between Conditions C1 --- C7 given in Subsection \ref{spoly}. In the proofs of Theorems \ref{mainTh}  and \ref{twopointsdesignth} we show that  these conditions are met under the conditions specified in the corresponding theorems. We need however first to work out the structure of the diagonal elements of the matrix $W_{11}$.

\subsection{The Diagonal of Matrix $W_{11}$}\label{Sw11}
 It will be convenient to use the following notation
\begin{align}
S_{1,n}&=\sum_{i=1}^n h_i,    \quad S_{2,n}=\sum_{i=1}^n h_i^2, \quad \text{and} \hfill \label{def4S1n} \\ 
V_n&=\frac{S_{2,n}}{n}-\left(\frac{S_{1,n}}{n}\right)^2=\frac{1}{n}\sum_{i=1}^n\left(h_i-\frac{S_{1,n}}{n}\right)^2 .\label{def4Vn}
\end{align}
Notice that if $S_{1,n}\neq0$ then \eqref{def4Vn} is equivalent to
\begin{align}
\frac{S_{2,n}}{S_{1,n}}-\frac{S_{1,n}}{n}&=\frac{V_n}{S_{1,n}/n}. \label{ABBn}
\intertext{Moreover, let us note that $d_{n+1}$ given by \eqref{def4d} reduces now to}
d_{n+1}&=\frac{1}{n}\left( \frac{\left(S_{1,n}- n h_{n+1}\right) ^{2}}{n^2 V_n }+1\right).\label{wzor4d}
\end{align}
Matrices $\AOOI $ and $W$ reduce  in the present case to the following form
\begin{align}
\AOOI &=\frac{1}{n^{2}V_n }\left[
\begin{array}{cc}
S_{2,n} & -S_{1,n} \\
-S_{1,n} & n%
\end{array}
\right], \label{def4A11lin}\\
W&=\gamma_W \left[
\begin{array}{cc}
\left( S_{2,n}-S_{1,n}h_{n+1}\right) ^{2} & -\left(
S_{2,n}-S_{1,n}h_{n+1}\right) \left( S_{1,n}-nh_{n+1}\right) \\
-\left( S_{2,n}-S_{1,n}h_{n+1}\right) \left( S_{1,n}-nh_{n+1}\right) &
\left( S_{1,n}-nh_{n+1}\right) ^{2}%
\end{array}
\right] ,\label{def4wlin}
\end{align}
where
\begin{align}
\gamma_W&=\frac{1}{q_{n+1}}\left(\frac{1}{n^{2}V_n }\right)^2.
\end{align}
and $q_{n+1}$ is given by \eqref{def4q}. We also get 
\begin{equation} \label{A004lin}
\AllI =\AOOI -W = \frac{1}{q_{n+1}}\frac{1}{n^{2}V_n }%
\left[
\begin{array}{cc}
S_{2,n+1} & -S_{1,n+1} \\
-S_{1,n+1} &  n+1
\end{array}%
\right]
\end{equation}
and
\begin{align}
\AO ^{T}D\AO &=D_{11}=\begin{bmatrix} \delta_{11} & \delta_{1h}\\
\delta_{1h} & \delta_{hh}
\end{bmatrix}
=\begin{bmatrix} \sum_{i=1}^n D_i & \sum_{i=1}^n D_ih_i\\
\sum_{i=1}^n D_ih_i & \sum_{i=1}^n D_ih_i^2
\end{bmatrix},\label{def4deltas}
\end{align}
where $D_{11}$ is given by \eqref{def4D11} and  $D_i$ are the diagonal elements of matrix $D$ given by \eqref{def4D} .
By applying \eqref{iidcase} we infer from \eqref{def4A11lin}, \eqref{def4wlin} and \eqref{A004lin} that
$q_{n+1}$ given by \eqref{def4q} equals
\begin{align}\label{qandV}
q_{n+1}&=\frac{(n+1)^2}{n^2}\frac{V_{n+1}}{V_n},
\end{align}
where $V_{n+1}$ stands for the {\em variance} of the design points given by \eqref{def4Vn} with substitution $n\leftarrow n+1$.
Clearly, since $D_i\ge 0$ the matrix  $\AO ^{T}D\AO $ is non-negative definite.
\begin{remark} It may be interesting to note that \eqref{qandV} combined with \eqref{def4q} and \eqref{wzor4d} is equivalent to the \ch{MR0143307} updating formula for $S_{2,n}$.
\end{remark}

It will be convenient to use the following notation.
\begin{align}\label{def4alphas}
\left(
\begin{array}{c}
a_{1}  \\
a_{2}       \\
a_{3}  \\
a_{4}       \\
a_{5}        \\
a_{6}
\end{array}
\right)=
\left(
\begin{array}{l}
q_{n+1}S_{2,n}  \\
S_{2,n+1}       \\
q_{n+1}S_{1,n}  \\
S_{1,n+1}       \\
nq_{n+1}        \\
n+1
\end{array}
\right)
& \quad \text{and}\quad
\left(
\begin{array}{c}
\alpha_1  \\
\alpha_2 \\
\alpha_3  \\
\alpha_4 \\
\alpha_5  \\
\end{array}
\right)
=
\left(
\begin{array}{l}
 a_1^2 -a_2^2 \\
 a_1a_3-a_2a_4\\
 a_3^2 -a_4^2 \\
 a_3a_5-a_4a_6\\
 a_5^2 -a_6^2 \\
\end{array}
\right).
\end{align}
We hasten to note that in the present case conditions \eqref{techCond4Wiipos}  are equivalent to \eqref{mainIneq1b} -\eqref{mainIneq2b}, respectively.  However, to get conditions equivalent to \eqref{techCond4Wiipos}  in a simpler algebraic form we begin with the following diagonal form of $W_{ii}$.

\begin{proposition} \label{diag4W11}
The diagonal elements of $W_{11}$ are given by
\begin{align}
W_{11}(1,1)&=\frac{\alpha_1\delta_{11}-2\alpha_2\delta_{1h}+\alpha_3\delta_{hh}}{(n+1)^4V_{n+1}^2}  \label{diag4W11_1}
\intertext{and}
W_{11}(2,2)&=\frac{\alpha_3\delta_{11}-2\alpha_4\delta_{1h}+\alpha_5\delta_{hh}}{(n+1)^4V_{n+1}^2}, \label{diag4W11_2}
\end{align}
respectively, where $\delta_{ij}$ are given by \eqref{def4deltas}. $\lozenge$
\end{proposition}

Prior to further considering conditions guaranteeing non-negativeness of the diagonals of $W_{11}$ we shall need some properties of polynomials driving the behavior of  numerators in \eqref{diag4W11_1} and \eqref{diag4W11_2}.

\subsection{Driving polynomials}\label{spoly}
Let us introduce two closely related pairs of quadratic polynomials: $\left(p_1(h),p_2(h)\right)$ and  $\left(\bar{p}_1(h),\bar{p}_2(h)\right)$, referred to in the following as driving polynomials. We need them to formulate conditions equivalent to the variance reduction property.
Polynomials $p_1(h)$ and $p_2(h)$ are given by
\begin{align}
p_1(h)&=\alpha_1 -2\alpha_2 h +\alpha_3 h^2 \quad \text{and} \label{p1}\\
p_2(h)&=\alpha_3 -2\alpha_4 h +\alpha_5 h^2, \label{p2}
\end{align}
respectively, where $\alpha_1, \alpha_2,\ldots,\alpha_5$ are given by \eqref{def4alphas}. In Lemmas \ref{two_deltas} and \ref{roots} we derive their roots
 $r_{1,1},r_{1,2}$ and $r_{2,1},r_{2,2}$ given by \eqref{r11}--\eqref{r22}, respectively and show features needed in the sequel.
 Polynomials  $\bar{p}_1(h)$ and $\bar{p}_2(h)$ are given by
\begin{align}
\bar{p}_1(h)&=\left(h-r_{1,1}\right)\left(h-r_{1,2}\right)\quad \text{and} \label{barp1}\\
\bar{p}_2(h)&=\left(h-r_{2,1}\right)\left(h-r_{2,2}\right) \label{barp2}
\end{align}
\noindent and are obtained from $p_1(h)$ and $p_2(h)$ by dividing them by $\alpha_3$ and $\alpha_5$, respectively.

\subsection{Necessary and sufficient conditions for the VRP. }\label{necsuff}
In the following we consider statements {\bf C1}-{\bf C7} and the corresponding assumptions under which they are equivalent to the non-negativeness of the diagonal of $W_{11}$.
\begin{description}
\item[C1] \label{C1} The diagonal elements of matrix $W_{11}$ given by \eqref{diag4W11_1}-\eqref{diag4W11_2} are non-negative for any vector of non-increasing and non-negative $D_i$'s.
\item[C2]  \label{C2} For all non-increasing and non-negative $D_i$'s
\begin{align}
\alpha_1\delta_{11}-2\alpha_2\delta_{1h}+\alpha_3\delta_{hh}&\ge 0, \quad \text{and}  \label{mainIneq1}\\
\alpha_3\delta_{11}-2\alpha_4\delta_{1h}+\alpha_5\delta_{hh}&\ge 0. \label{mainIneq2}
\end{align}
\item[C3]  \label{C3} For all non-increasing and non-negative $D_i$'s
\begin{align}
\sum_{i=1}^n D_i\left( \alpha_1 -2\alpha_2h_i+\alpha_3h_i^2 \right)&\ge 0 \quad \text{and} \label{mainIneq1a}\\
\sum_{i=1}^n D_i\left( \alpha_3-2\alpha_4h_i+\alpha_5h_i^2\right) &\ge 0. \label{mainIneq2a}
\end{align}
\item[C4]  \label{C4} For each $m\in\left\{1,...,n\right\}$
 \begin{align}
  \alpha_1 -2\alpha_2\frac{S_{1,m}}{m}+\alpha_3\frac{S_{2,m}}{m} &\ge 0, \quad \text{and} \label{mainIneq1b}\\
\alpha_3 -2\alpha_4\frac{S_{1,m}}{m}+\alpha_5\frac{S_{2,m}}{m} &\ge 0.  \label{mainIneq2b}
 \end{align}
\item[C5] \label{C5} For each $m\in\left\{1,...,n\right\}$  we have
\begin{align}
  p_1\left(\frac{S_{1,m}}{m}\right)+\alpha_3 V_m  &\ge 0, \quad \text{and} \label{mainIneq1c}\\
 p_2\left(\frac{S_{1,m}}{m}\right)+\alpha_5 V_m  &\ge 0, \label{mainIneq2c}
\end{align}
where polynomials $p_1(h)$ and $p_2(h)$ are given by \eqref{p1} and \eqref{p2}, respectively.
\item[C6]  \label{C6} For each $m\in\left\{1,...,n\right\}$  we have
    \begin{align}
\bar{p}_1\left(\frac{S_{1,m}}{m}\right)+ V_m  &\ge 0, \quad \text{and} \label{mainIneq1d}\\
\bar{p}_2\left(\frac{S_{1,m}}{m}\right)+ V_m  &\ge 0, \label{mainIneq2d}
\end{align}
where polynomials $\bar{p}_1(h)$ and $\bar{p}_2(h)$ are given by \eqref{barp1} and \eqref{barp2}, respectively.
\item[C7]  \label{C7} For each $m\in\left\{1,...,n\right\}$
\begin{align}
\bar{p}_1\left(\frac{S_{1,m}}{m}\right)+ V_m  &\ge 0. \label{mainIneq1e}
\end{align}
\end{description}
\begin{theorem}\label{equivCond}
\begin{enumerate}
\item \label{part1} If $V_n>0$ then statements {\bf C1 -- C5} imply each other.
\item \label{part2} If $V_n,\alpha_3$ and $\alpha_5$ are positive then {\bf C1 -- C6} imply each other.
\item \label{part3} If $h_1,\ldots,h_{n+1}$ are non-negative and increasing then {\bf C1 -- C7} imply each other.
\end{enumerate}
\end{theorem}
\ofproof{Theorem}{\ref{equivCond}}{Assuming that $V_n>0$, the equivalence of {\bf C1, C2} and {\bf C3} is evident. Lemma \ref{Pos}
applied with $U_i=\alpha_1 -2\alpha_2h_i+\alpha_3h_i^2 $ implies the equivalence of {\bf C3} and {\bf C4}. By applying \eqref{def4Vn} we
get equivalence of {\bf C4} and {\bf C5}. Assuming additionally that $\alpha_3>0$ and $\alpha_5>0$ we can divide \eqref{mainIneq1c} by $\alpha_3$
and \eqref{mainIneq2c} by $\alpha_5$, respectively, without changing the direction of the inequalities. This shows equivalence of {\bf C5} and {\bf C6}.
If  $h_1,\ldots,h_{n+1}$ are non-negative and increasing then  by \eqref{r11equalsr21}--\eqref{r12-r22}, Lemma \ref{scopehn1} and Corollaries \ref{rootsup3} and \ref{rootsup4} both polynomials have the same roots $r_{1,1}=r_{2,1}$. Moreover, the  root $r_{1,2}$ of $\bar{p}_1$ is larger than the root $r_{2,2}$ of $\bar{p}_2$ and $r_{2,2}\ge r_{2,1}$. This implies that $\bar{p}_2(h)\ge \bar{p}_1(h)$ for $h \in [r_{1,1},r_{1,2}]$. Hence for $h \in [r_{1,1},r_{1,2}]$ \eqref{mainIneq1d} implies \eqref{mainIneq2d}. Since both $\bar{p}_1(h)$ and $\bar{p}_2(h)$ are non-negative for $h \notin [r_{1,1},r_{1,2}]$ we get the equivalence of {\bf C6} and {\bf C7}. This concludes the proof of Theorem \ref{equivCond}. }

\begin{remark} \label{r2}
Notice that the definition of $\alpha_i$'s given in \eqref{def4alphas} implies that the left hand sides of inequalities \eqref{mainIneq1b}--\eqref{mainIneq2b} in statement {\bf C4} are polynomials of degree 4 of variable $h_{n+1}$. Hence it is clear that, given $h_1,\ldots,h_n$, the set of $h_{n+1}$'s leading to the variance reduction equals the intersection of the positivity regions of $2m$ polynomials of degree 4. No transparent characterization of these sets seems available, yet for numerical values of $h_1,\ldots,h_n$ it is always possible to determine, at least numerically, which values of $h_{n+1}$ need to be avoided to retain the VRP of the design. This opens a way to practical sequential methods of determining designs with the VRP for decreasing variances.
\end{remark}

\subsection{Proofs o Theorems \ref{mainTh} and \ref{twopointsdesignth}.}


\ofproof{Theorem}{\ref{mainTh}}{By Theorem \ref{equivCond} it is enough to show that statement {\bf C7} holds true.

Clearly \eqref{mainIneq1e} holds true for $m$ such that $\frac{S_{1,m}}{m}<r_{1,1}$. Let $m_1$ be the smallest $m$ such that
$\frac{S_{1,m}}{m}\ge r_{1,1}$. Since $h_i$ are increasing $\frac{S_{1,m}}{m}$ is increasing with $m$ and hence we have
\begin{align*}
\frac{S_{1,m}}{m}\ge r_{1,1} \quad \text{for $m \in [m_1,n]$}.
\end{align*}
Inequalities \eqref{Snbetweenroots} imply that the interval $[m_1,n]$ is non-empty, it contains at least the right-hand end point $n$.

Let us note that
\begin{align*}
\bar{p}_1\left(\frac{S_{1,m}}{m}\right)+ V_m &= \left(\frac{S_{1,m}}{m}-r_{1,1}\right)\left(\left(\frac{S_{1,m}}{m}-r_{1,1}\right)+\frac{V_m}{\left(\frac{S_{1,m}}{m}-r_{1,1}\right)}-\left(r_{1,2}-r_{1,1}\right)\right). \label{wazne}
\end{align*}
From Lemmas  \ref{lbarp} and \ref{scopehn1}  we infer that
\begin{equation*}
\left(\frac{S_{1,m}}{m}-r_{1,1}\right)+\frac{V_m}{\left(\frac{S_{1,m}}{m}-r_{1,1}\right)}>\left(r_{1,2}-r_{1,1}\right) \label{wazne2}
\end{equation*}
holds true for $m=n$. We will show that for $m \in [m_1,n]$ the left hand side of \eqref{wazne2} is decreasing with $m$ or, equivalently, that
the expression on the right hand side of \eqref{form4goingdown} is positive. Hence, we need to show that
\begin{align}
h_{m+1}<\frac{S_{2,m}-r_{1,1}S_{1,m}}{S_{1,m}-mr_{1,1}}.\label{wazne4}
\end{align}
Let us fix $m\ge m_1$ and define function $g(r)$ by
\begin{align}
 g(r) &= \frac{S_{2,m}-rS_{1,m}}{S_{1,m}-mr}.\label{def4g}
 \end{align}
By Corollary \ref{rootsup3} we have
\begin{align}
 r_{1,1}\left( m\right) =g(h_{m+1}) \le r_{1,1} < \frac{S_{1,m}}{m}. \label{wazne3}
\end{align}
Notice that $h_{m+1}=g\left(r_{1,1}\left( m\right)\right)$. Since for $r< \frac{S_{1,m}}{m}$ the function $g(r)$ is increasing  \eqref{wazne3} implies \eqref{wazne4}. Hence we get that
\[ \bar{p}_1\left(\frac{S_{1,m}}{m}\right)+ V_m > 0\]
for any $m \in [m_1,n]$. This concludes our proof of  Theorem \ref{mainTh}. }

\ofproof{Theorem}{\ref{twopointsdesignth}}{
The design points $h_i$ in the present case, in contrast with the case considered in Theorem \ref{mainTh}, need be neither non-negative nor increasing and $\alpha_3$ can be negative, eg for $b=-1,c=1$. By part \ref{part1} of Theorem \ref{equivCond} we can rely  only on the equivalence of statements   {\bf C1 -- C5}.
To show that this two-point regression design has the VRP  we need to consider two cases: for $n$ even and for $n$ odd.
\begin{description}
\item[The case of even $n$.]
Condition \eqref{mainIneq1c} reduces in the present case to
\begin{align*}
p_{1}\left( \frac{S_{1,m}}{m}\right) +\alpha _{3}V_{m}&=\left\{
\begin{array}{ccc}
 \frac{1}{2}\left( n+1\right)c^{2}\left( b-c\right)^{2}  & \text{\ if\ } & m\text{\ is even} \\
\frac{1}{2}\frac{m+1}{m}\left( n+1\right) c^{2}\left( b-c\right) ^{2}  & \text{\ if\ } & m\text{\ is odd}%
\end{array}%
\right. \nonumber \\
& >0 \quad \text{for} \quad m\in \left[ 1,n\right] .
\end{align*}

Condition \eqref{mainIneq2c} reduces in the present case to
\begin{align*}
p_{2}\left( \frac{S_{1,m}}{m}\right) +\alpha _{5}V_{m}&=\left\{
\begin{array}{ccc}
\frac{1}{2}\left( n+1\right)\left( b-c\right) ^{2} & \text{\ if\ } & m\text{\ is even} \\
\frac{1}{2}\frac{m+1}{m}\left(n+1\right) \left( b-c\right) ^{2}  & \text{\ if\ } & m\text{\ is odd}%
\end{array}
\right. 
\end{align*}
and clearly these expressions in are positive for distinct $b$ and $c$ and $m\in \left[ 1,n\right]$.

\item[The case of odd $n$.]
Conditions \eqref{mainIneq1c}-\eqref{mainIneq2c} reduce in the present case to
\begin{align*}
p_{1}\left( \frac{S_{1,m}}{m}\right) +\alpha _{3}V_{m}&=\left\{
\begin{array}{ccc}
\frac{1}{2}n\frac{\left( n+1\right) ^{2}}{\left( n-1\right)
^{2}}b^{2}\left( b-c\right) ^{2} & \text{\ if\ } & m\text{\ is even}
\\
\frac{1}{2}\frac{m-1}{m} n\frac{\left(
n+1\right) ^{2}}{\left( n-1\right) ^{2}}b^{2}\left( b-c\right) ^{2} & \text{\ if\ } & m%
\text{\ is odd}%
\end{array}%
\right. \nonumber \\
&\geq 0 \quad \text{for} \quad m\in \left[ 1,n\right] 
\intertext{and similarly}
p_{2}\left( \frac{S_{1,m}}{m}\right) +\alpha _{5}V_{m}&=\left\{
\begin{array}{ccc}
\frac{1}{2}n\frac{\left( n+1\right) ^{2}}{\left(
n-1\right) ^{2}}\left( b-c\right) ^{2} & \text{\ if\ } & m~%
\text{\ is even} \\
\frac{1}{2}\frac{m-1}{m}%
n\frac{\left( n+1\right) ^{2}}{\left( n-1\right)
^{2}}\left( b-c\right) ^{2} & \text{\ if\ } & m~%
\text{\ is odd}
\end{array}%
\right.\nonumber \\
&\geq 0 \quad \text{for} \quad m\in \left[ 1,n\right].
\end{align*}
\end{description}

So, both for $n$ even and for $n$ odd the condition {\bf C5} is met and hence by Part \ref{part1} of Theorem \ref{equivCond} we conclude that the considered two-point design has the VRP.}

\section{Examples and counterexamples}\label{examples}
\begin{example} \label{ex01} \sl Matrix $W_{11}$ need not be positive definite, even for decreasing $\sigma^2_i$ and increasing $h_i$'s. Let
\begin{align*}
\Al =\left[
  \begin{array}{cc}
    1 & 0.62 \\
    1 & 1.24 \\
    1 & 1.80 \\
    1 & 1.96 \\
  \end{array}
\right]
\quad \text{and}& \quad
\Sgmall = Diag\left\{1.56,1.26,0.78,0.28\right\}.
\end{align*}
Then we have
\begin{align*}  
W_{11} = \left[
           \begin{array}{cc}
            0.54605 & -0.55859 \\
             -0.55859 & 0.56787 \\
           \end{array}
         \right]
\end{align*}
and the eigenvalues of $W_{11}$ are given by
\[  -0.00174,\quad 1.11566 ,\]
respectively. So, the matrix $W_{11}$ is indefinite, yet the diagonal elements are positive.
\end{example}

\begin{example} \label{ex02} \sl For correlated $Y_1,\ldots,Y_n$ and uncorrelated $Y_{n+1}$ the VRP need not hold  even for decreasing $\sigma^2_i$ and increasing $h_i$'s. Let
\begin{align*}
\Al =\left[
\begin{array}{cc}
1 & \frac{7}{10} \\
1 & \frac{8}{5} \\
1 & \frac{17}{10} \\
1 & \frac{19}{10}%
\end{array}%
\right],
\quad
S&=\left[
\begin{array}{ccc}
2 & \frac{5}{4} & \frac{9}{10} \\
\frac{5}{3} & \frac{6}{5} & \frac{2}{5} \\
\frac{5}{4} & 1 & \frac{1}{4}%
\end{array}%
\right],
\quad
\Sgmall =\left[
\begin{array}{cc}
S^{T}S & \mathbf{0} \\
\mathbf{0}' & \frac{1}{5}
\end{array}
\right] ,
\end{align*}
where $\mathbf{0}'=[0,0,0]$. Then the diagonal of $\Sgmall$ equals $\left\{8.340\,278,4.0025,1.0325,0.2\right\}$ and
\begin{align*}
V_n-V_{n+1}&=  \left[
\begin{array}{cc}
-4. 505\,055 &
3. 277\,313 \\
3. 277\,313 &
-1. 987\,787
\end{array}
\right],
\end{align*}
ie. adding the fourth observation results in an increase of variance for both intercept and slope estimators. We refer to Proposition \ref{t2dep} for a general decomposition in the case of correlated random variables.
\end{example}

\begin{example} \label{ex03} \sl For uncorrelated $Y_1,\ldots,Y_n,Y_{n+1}$ the VRP need not hold  even for decreasing $\sigma^2_i$ and increasing $h_i$'s $i=1,\ldots,n$ but with $h_{n+1}<h_n$. Let
\begin{align*}
\Al =\left[
\begin{array}{cc}
1 & 0.7 \\
1 & 1.6 \\
1 & 1.62 \\
1 & 1.45
\end{array}
\right]
\quad \text{and} \quad
\Sgmall =\left[
\begin{array}{cccc}
2 & 0 & 0 & 0 \\
0 & 1 & 0 & 0 \\
0 & 0 & \frac{4}{5} & 0 \\
0 & 0 & 0 & \frac{1}{5}%
\end{array}%
\right] .
\end{align*}
Then we have
\begin{align*}
V_n - V_{n+1} &=\left[
\begin{array}{cc}
-0.002\,878  &
0.064\,592 \\
0.064\,592 &
-0.013\,034
\end{array}
\right].
\end{align*}
\end{example}

\appendix
\section{Auxiliary results.} \label{app}
In derivation of our results we use the inversion formulae for block partitioned square matrices, included below for the convenience of reader:
\begin{itemize}
\item assuming that the inverse matrices exist (cf. \ch{MR0346957}, Problem 1.2.8 ), {\small
\begin{equation}\label{p128}
\hspace{-0.4in}\left[
\begin{array}{cc}
C & B \\
B^{T } & D%
\end{array}%
\right] ^{-1}=\left[
\begin{array}{cc}
C^{-1}+C^{-1}B\left( D-B^{T }C^{-1}B\right) ^{-1}B^{T }C^{-1} &
-C^{-1}B\left( D-B^{T }C^{-1}B\right) ^{-1} \\
-\left( D-B^{T }C^{-1}B\right) ^{-1}B^{T }C^{-1} & \left(
D-B^{T }C^{-1}B\right) ^{-1}%
\end{array}%
\right]
\end{equation}}
\item  in the case of a nonsingular matrix $C$ and two column-vectors $B$ and $D$ (cf. \ch{MR0346957}, Problem 1.2.7)
\begin{equation}\label{p127}
\hspace{-1in}\left(C + BD^T\right)^{-1}= C^{-1}- \frac{1}{1+D^TCB}\left(C^{-1}B\right)\left(D^TC^{-1}\right).
\end{equation}
\end{itemize}

{\bf Proof of Proposition \ref{t2}. } Let us note that by applying \eqref{p127} and partitions \eqref{m2c} one gets from \eqref{def4V00} the well known \ch{MR0036980} updating formula \eqref{iidcase} for covariances of the LS estimators in the iid case, which we use in the sequel.
Let us also note that
\begin{align}
W\AO ^{T}\AO W &= \frac{a\AOOI a^{T}}{1+a\AOOI a^{T}}W, &&
&Wa^{T}aW &= \frac{\left(a\AOOI a^{T}\right) ^{2}}{1+a\AOOI a^{T}}W, \label{prosto1}\\
Wa^{T}a\AOOI &= a\AOOI a^{T}W &\text{and}&
&\AOOI a^{T}aW &= a\AOOI a^{T}W.\label{prosto2}
\end{align}
By applying next the block-matrix inversion formula \eqref{p128} and using  \eqref{def4D} one can derive  \eqref{kwadratyporedukcjidoI} with some algebra, omitted for the sake of compactness. $\lozenge$\\

The following simple lemma is pivotal for the present paper.
\begin{lemma} \label{Pos} Let $U_i, i=1,2,\ldots,n$ be given.  Inequalities
\begin{align}
\sum_{i=1}^n D_iU_i &\ge 0 \label{PosD}
\end{align}
hold true for all non-negative and non-increasing $D_i$'s if and only if
\begin{align}
\sum_{i=1}^m U_i &\ge 0 \quad \text{for $m=1,\ldots,n$}. \label{Posk}
\end{align}
\end{lemma}
\ofproof{Lemma}{\ref{Pos}}{Clearly, by choosing $D_i=1$ for $i=1,\ldots,m$ and $D_i=0$ for $i>m$ we find that \eqref{PosD} implies \eqref{Posk}.
To show that \eqref{Posk}  implies \eqref{PosD} let us note that
\begin{align} \label{PosM}
\sum_{i=1}^n D_iU_i &= \sum_{j=1}^{n-1}\left( D_{j}-D_{j+1}\right)\sum_{s=1}^{j}U_{s}+D_{n}\sum_{s=1}^{n}U_{s}
\end{align}
and hence that the right hand side of \eqref{PosM} is non-negative if \eqref{Posk} holds true. }
\begin{remark} Property \eqref{Posk} should not be confused with the majorization of \ch{MR0046395} because no ordering of $U_i$'s is here assumed.
\end{remark}

\noindent {\bf Proof of Proposition \ref{char4WiiPos}.} By \eqref{def4W11} we have
\begin{align*}
W_{ii} &= \sum_{j=1}^n \left(\sum_{l=1}^k \AOOI (i,l)\AO (j,l)\right)^2 D_j^2 -\sum_{j=1}^n \left(\sum_{l=1}^k \AllI (i,l)\AO (j,l)\right)^2 D_j^2\\
&= \sum_{j=1}^n \left( \left(\sum_{l=1}^k \AOOI (i,l)\AO (j,l)\right)^2  -\sum_{j=1}^n \left(\sum_{l=1}^k \AllI (i,l)\AO (j,l)\right)^2 \right) D_j^2 \ge0
\end{align*}
Hence, in order to have
\[ W_{ii} \ge 0 \]
for every non-increasing $D_j^2$ with   $D_j^2=\sigma_i^2 - \sigma_{n+1}^2$, it is sufficient and necessary, by Lemma \ref{Pos}, to have \eqref{techCond4Wiipos}. This completes the proof.~$\lozenge$\\

\begin{lemma} \label{two_deltas}
We have
\begin{align}
\Delta_1 &=  \alpha_2^2 - \alpha_1\alpha_3 \nonumber\\
&=  q_{n+1}^{2}\left( S_{2,n}-S_{1,n}h_{n+1}\right)
^{2}h_{n+1}^{2}\ge 0 \label{Delta2}
\intertext{and}
\Delta_2 &= \alpha_4^2 -\alpha_3\alpha_5 \nonumber \\
&= q_{n+1}^{2}\left( S_{1,n}-nh_{n+1}\right) ^{2} \ge 0. \label{Delta3}
\end{align}
\end{lemma}
{\bf Proof.} The proof of Lemma \ref{two_deltas} follows immediately from \eqref{def4alphas} and from the identity 
\[\left(\alpha \beta-\delta\gamma\right)^2-\left(\alpha^2-\delta^2\right)\left(\beta^2-\gamma^2\right)=%
\left(\alpha\gamma-\beta\delta\right)^2.\quad \lozenge\]
Hence, Lemma \ref{two_deltas} implies that both polynomials $p_1(h)$ and $p_2(h)$ have roots. We will always assume that denominators differ from zero.

\begin{lemma}\label{roots}
Polynomials \eqref{p1} and \eqref{p2} have roots
\begin{align}
r_{1,1}&=\frac{\alpha _{2}}{\alpha _{3}}-\frac{\sqrt{\Delta _{1}}}{\alpha%
_{3}}=\frac{a_{1}-a_{2}}{%
a_{3}-a_{4}}=\frac{q_{n+1}S_{2,n}-S_{2,n+1}}{q_{n+1}S_{1,n}-S_{1,n+1}}, \label{r11}\\
r_{1,2}&=\frac{\alpha _{2}}{\alpha _{3}}+\frac{\sqrt{\Delta _{1}}%
}{\alpha _{3}}=\frac{a_{1}+a_{2}}{%
a_{3}+a_{4}}=\frac{q_{n+1}S_{2,n}+S_{2,n+1}}{q_{n+1}S_{1,n}+S_{1,n+1}} , \label{r12}
\intertext{and}
r_{2,1}&=\frac{\alpha _{4}}{\alpha _{5}}-\frac{\sqrt{\Delta_2}}{\alpha _{5}}=\frac{a_{3}-a_{4}}{a_{5}-a_{6}}=\allowbreak \frac{%
S_{1,n+1}-q_{n+1}S_{1,n}}{n+1-nq_{n+1}}, \label{r21}\\
r_{2,2}&=\frac{\alpha _{4}}{\alpha _{5}}+\frac{\sqrt{\Delta_2}}{\alpha _{5}}=\frac{a_{3}+a_{4}}{a_{5}+a_{6}}=\allowbreak \frac{%
q_{n+1}S_{1,n}+S_{1,n+1}}{n+1+nq_{n+1}}, \label{r22}
\end{align}
where $\Delta_1$ and $\Delta_2$ are given by \eqref{Delta2}  and \eqref{Delta3}, respectively.  Moreover, we have
\begin{align}
r_{1,1}=r_{2,1}&=\frac{S_{2,n}-S_{1,n}h_{n+1}}{S_{1,n}-nh_{n+1}} \label{r11equalsr21},\\
r_{1,2}-r_{1,1}&=2nV_n\frac{h_{n+1}}{\frac{nV_n}{S_{1,n}}+\left( h_{n+1}-\frac{S_{2,n}}{S_{1,n}}\right) }%
\frac{\left( n+1\right)n^{2} V_n +n^{2}\left(%
h_{n+1}-\frac{S_{1,n}}{n}\right) ^{2}}{n^{2} V_n \left( S_{1,n}+2nS_{1,n}+nh_{n+1}\right) +n^{2}S_{1,n} \left( h_{n+1}-%
\frac{S_{1,n}}{n}\right) ^{2}}, \label{r12-r11}
\intertext{and}
r_{1,2}-r_{2,2}&=\frac{2V_{n}\left( n\left( h_{n+1}-\frac{S_{1,n}}{n}\right) ^{2}+\left(
n^{2}+1\right) V_{n}\right) }{S_{1,n}\left( \left( h_{n+1}-\frac{S_{1,n}}{n}%
\right) ^{2}+h_{n+1}\frac{V_{n}}{nS_{1,n}}+\left( 2n+1\right) V_{n}\right) }%
\frac{\left( h_{n+1}-\frac{S_{1,n}}{n}\right) ^{2}+\left( 2n+1\right) V_{n}}{%
\left( \left( h_{n+1}-\frac{S_{1,n}}{n}\right) ^{2}+2\left( n+1\right)
V_{n}\right) }.\label{r12-r22}
\end{align}
\end{lemma}
Let us note that Lemmas \ref{two_deltas}--\ref{rinc1} are valid for any $h_1,\ldots,h_{n+1}$ for which the corresponding denominators differ from zero and they have been obtained by purely algebraic manipulations. For increasing non-negative sequences $h_1,\ldots,h_{n+1}$ Lemma \ref{scopehn1} implies that the right hand sides of \eqref{r11equalsr21}-\eqref{r12-r11} are positive.

 We shall need the following lemmas.
\begin{lemma}\label{lbarp}
We have
\begin{align}
\bar{p}_1\left(\frac{S_{1,n}}{n}\right)+ V_n  &= nV_{n}\frac{\left( S_{2,n}-S_{1,n}h_{n+1}\right) }{\left( S_{1,n}-h_{n+1}n\right) }\frac{%
\left( h_{n+1}-\frac{S_{1,n}}{n}\right) ^{2}+\left( 2n+1\right) V_{n}}{S_{1,n}\left( h_{n+1}-\frac{S_{1,n}%
}{n}\right) ^{2}+\left( nh_{n+1}+S_{1,n}\left( 2n+1\right) \right) V_{n}} \label{mainIneq1f}
\end{align}
\end{lemma}

\begin{lemma}\label{goingdown}
For $m < n$ we have
\begin{multline}
\left( \left( \frac{S_{1,m}}{m}-r_{1,1}\right) +\frac{V_m }{\left( \frac{S_{1,m}}{m}-r_{1,1}\right) }\right) -\left(%
\left( \frac{S_{1,m+1}}
{m+1}-r_{1,1}\right) +\frac{V_{m+1}}{\left( \frac{S_{1,m+1}}{m+1}-r_{1,1}\right) }\right) \\
= -\frac{ h_{m+1}-r_{1,1} }{%
\left( m+1\right) \left( \frac{S_{1,m+1}}{m+1}-r_{1,1}\right) }\left( h_{m+1}-\frac{S_{2,m}-r_{1,1}S_{1,m}}{S_{1,m}-mr_{1,1}}%
\right) . \label{form4goingdown}
\end{multline}
\end{lemma}

\begin{lemma} \label{rinc1}
For $m < n$ we have
\begin{align}
&\frac{S_{2,m+1}-\left( S_{1,m+1}\right) h_{m+2}}{S_{1,m+1}-\left(
m+1\right) h_{m+2}}-\frac{S_{2,m}-S_{1,m}h_{m+1}}{S_{1,m}-mh_{m+1}} \nonumber\\
&=\left(h_{m+2}- h_{m+1}\right) \frac{m\left(
\frac{S_{1,m}}{m}-h_{m+1}\right) ^{2}+m\left( m+1\right) V_{m}}{\left( S_{1,m}-h_{m+1}m\right)
\left(S_{1,m+1}-(m+1)h_{m+2}\right) }.\label{rup}
\end{align}
\end{lemma}

We shall also need the following properties of the design points $h_i$ which can be easily derived using convexity arguments.

\begin{lemma}\label{scopehn1}
If $n>1$ and  $h_i$ are non-negative and increasing with $i$ then
 $\alpha_1,\ldots,\alpha_5$ are positive,
\begin{align}
 h_1 < \frac{S_{1,n}}{n} & <\frac{S_{2,n}}{S_{1,n}}<h_{n}, \label{preincreaseh}\\
 h_1 < r_{1,1}  < \frac{S_{1,n}}{n} & < r_{2,2}<\frac{S_{1,n+1}}{n+1} < h_{n+1}, \label{Snbetweenroots}
 \intertext{and}
\frac{S_{2,n}}{S_{1,n}} & < r_{1,2}<\frac{S_{2,n+1}}{S_{1,n+1}} < h_{n+1}. 
\end{align}
\end{lemma}
\ofproof{lemma}{\ref{scopehn1}}{Since
\[ \frac{S_{2,n}}{S_{1,n}}=\sum_{i=1}^n\frac{h_i}{S_{1,n}}h_i \]
we get
\[\frac{S_{1,n}}{n}< \frac{S_{2,n}}{S_{1,n}}<\max\{h_i\}=h_n.\]
Indeed, the left hand side inequality holds because heavier weights are assigned to $h_i$'s of higher value. This implies \eqref{preincreaseh}.
Moreover, we have
\begin{align*}
r_{1,1}&=\frac{%
S_{1,n}h_{n+1}-S_{2,n}}{nh_{n+1}-S_{1,n}}
=\sum_{i=1}^{n}\frac{\left( h_{n+1}-h_{i}\right) }{\sum_{j=1}^{n}\left(
h_{n+1}-h_{j}\right) }h_{i}
\intertext{and hence}
h_1&<r_{1,1}<\frac{S_{1,n}}{n}
\intertext{and the right hand side inequality holds because heavier weights are assigned to $h_i$'s of smaller values. Clearly, we have $h_1<r_{1,1}$. Let us note that}
r_{1,2}&=\frac{q_{n+1}S_{2,n}+S_{2,n+1}}{q_{n+1}S_{1,n}+S_{1,n+1}}=\frac{q_{n+1}S_{1,n}}{q_{n+1}S_{1,n}+S_{1,n+1}}%
\frac{S_{2,n}}{S_{1,n}}+\frac{S_{1,n+1}}{q_{n+1}S_{1,n}+S_{1,n+1}}\frac{%
S_{2,n+1}}{S_{1,n+1}}.
\intertext{Hence }
\frac{S_{2,n}}{S_{1,n}} &< r_{1,2}<\frac{S_{2,n+1}}{S_{1,n+1}}.
\intertext{Finally, we have}
r_{2,2}&=\frac{q_{n+1}S_{1,n}+S_{1,n+1}}{n+1+nq_{n+1}}=\frac{nq_{n+1}}{n+1+nq_{n+1}}%
\frac{S_{1,n}}{n}+\frac{n+1}{n+1+nq_{n+1}}\frac{S_{1,n+1}}{n+1}
\end{align*}
and hence
\[ \frac{S_{1,n}}{n} < r_{2,2} < \frac{S_{1,n+1}}{n+1}.\]
}

Lemmas \ref{rinc1} and \ref{scopehn1} imply the following.
\begin{corollary}\label{rootsup3}
Let
\begin{align}
r_{1,1}\left(m\right) &= \frac{S_{2,m}-S_{1,m}h_{m+1}}{S_{1,m}-mh_{m+1}}=g(h_{m+1}),
\end{align}
where $g$ is given by \eqref{def4g}. If $h_i$ are non-negative and increasing then for  $m<n$ then we have
\begin{equation*}
0 \le r_{1,1}\left( m\right) \le r_{1,1}\left(m+1\right)  \le r_{1,1}\left(n\right)=r_{1,1}.
\end{equation*}
\end{corollary}

\begin{corollary}\label{rootsup4}
For
$h_1,\ldots,h_{n+1}$ non-negative and increasing the values of all expressions on the right hand sides at  \eqref{r11equalsr21} ---  \eqref{rup} are strictly positive.
\end{corollary}
Let us finally note that in the case of correlated random variables we have the following decomposition extending \eqref{kwadratyporedukcjidoI} onto the case of general correlated random variables. We keep here the same notation
as in Proposition \ref{t2}.

\begin{proposition}
\label{t2dep} If the covariance matrix $\mathbf{\Sgmall }$ and design matrix $\Al $ in a model given by \eqref{m2a} and \eqref{m2c}  are of full ranks  then we have
\begin{align}
\Vll  =
\VOO -\mathbf{\sigma }_{n+1}^{2}W-W_{11}+W_{22}, \label{Decomp4GenVarDep}
\end{align}
where
\begin{align}
W_{22} &= \AllI \left( a^{T}\mathbf{\bar{\sigma}}_{1,n+1}^{\prime}\AO +\AO ^{T}\mathbf{\bar{\sigma}}_{1,n+1}a\right) \AllI ,
\end{align}
and $\bar{\sigma}_{1,n+1}$ is the $(n,1)$ vector of covariances of $Y_{n+1}$ and $\left(Y_1,Y_2,...,Y_n\right)$.
\end{proposition}

Proposition \ref{t2dep} can be shown in a similar way like the decomposition \eqref{kwadratyporedukcjidoI} of Proposition \ref{t2}.   Proposition \ref{t2dep} shows explicitly how the difference between $\Vll $ and $\VOO $ is affected by covariances. We also refer to Example \ref{ex02}, where consequences of decomposition \eqref{Decomp4GenVarDep} are shown and where the non-zero correlations affect the VRP. 

\bibliography{VR4OLS0}
\end{document}